\documentclass[10pt]{article}

\newcounter{mascotsection}
\newcounter{mascotsubsection}[mascotsection]

\newcommand{\mascotsection}[1]{%
	\stepcounter{mascotsection}
	\section{\themascotsection . #1}
}
\newcommand{\mascotsubsection}[1]{%
	\stepcounter{mascotsubsection}
	\subsection{\themascotsection .\themascotsubsection . #1}
}

\usepackage{epsfig}

\usepackage{amsfonts,amsmath,amssymb,amsthm}
\usepackage{graphicx}

\newcommand{\NN}{\mathbb{N}}
\newcommand{\RR}{\mathbb{R}}

\begin{document}
\vskip 1.8cm 
\centerline{\LARGE\bf Two-dimensional interpolation using}
\vskip 0.3cm
\centerline{\LARGE\bf a cell-based searching procedure}
\vskip 0.6cm
\centerline{\bf Roberto Cavoretto}
\centerline{Department of Mathematics \lq\lq G. Peano\rq\rq, University of Torino,}
\centerline{Via Carlo Alberto 10,~10123 Torino,~Italy}
\centerline{\em roberto.cavoretto@unito.it}
\vskip .6cm
\begin{abstract}
\noindent In this paper we present an efficient algorithm for bivariate interpolation, which is based on the use of the partition of unity method for constructing a global interpolant. It is obtained by combining local radial basis function interpolants with locally supported weight functions. In particular, this interpolation scheme is characterized by the construction of a suitable partition of the domain in cells so that the cell structure strictly depends on the dimension of its subdomains. This fact allows us to construct an efficient cell-based searching procedure, which provides a significant reduction of CPU times. Complexity analysis and numerical results show such improvements on the algorithm performances.
\end{abstract}

\medskip

\medskip
\mascotsection{Introduction}
Let $\{(\textbf{x}_i,f_i), i=1,2,\ldots,n\}$ be a finite set of discrete data, with $\textbf{x}_i\in \Omega \subseteq \RR^2$, and $f_i \in \RR$. The $\textbf{x}_i$ are called the \textsl{nodes}, while the (corresponding) $f_i$ are the \textsl{data values}. The latter are obtained by sampling some (unknown) function $f:\Omega \rightarrow \RR$ at the nodes, i.e. $f_i=f(\textbf{x}_i)$, $i=1,2,\ldots,n$. 

Therefore, the scattered data interpolation problem consists in finding a continuous function ${\cal R}:\Omega \rightarrow \RR$ such that
\begin{equation}  \label{intcond}
{\cal R}(\textbf{x}_i)=f_i, \hspace{1cm} i=1,2,\ldots,n.
\end{equation}
Here, we consider the problem of constructing an efficient algorithm for bivariate interpolation of (large) scattered data sets. It is based on the \textsl{partition of unity method} for constructing a global interpolant by blending radial basis functions (RBFs) as local approximants and using locally supported weight functions.
 
Now, starting from the results of previous researches (see \cite{Allasia11,Cavoretto10a,Cavoretto12a,Cavoretto12b, Cavoretto12c}) where efficient searching procedures based on the partition of the domain in strips or spherical zones are considered, we extend the previous ideas replacing the strip-based partition structure with a cell-based one. The latter leads to the creation of a cell-based searching procedure, whose origin comes from the repeated use of a quicksort routine with respect to different directions, enabling us to pass from not ordered to ordered data structures. In particular, this process is strictly related to the construction of a partition of the domain $\Omega$ in square cells, which consists in generating two orthogonal families of parallel strips, where the original data set is suitably split up in ordered data subsets.

Then, exploiting the ordered data structure and the domain partition, the cell algorithm is efficiently implemented and optimized by connecting the me\-thod itself with the effective cell-based searching procedure. More precisely, the considered technique is characterized by the construction of a double structure of crossed strips, which partitions the domain in square cells and strictly depends on the dimension of its subdomains, providing a significant improvement in the searching procedures of the nearest neighbour points compared to the searching techniques in \cite{Allasia11, Cavoretto10a, Cavoretto12a}. The final result is an efficient algorithm for bivariate interpolation of scattered data. Finally, complexity analysis and numerical tests show the high efficiency of the proposed algorithm.

The paper is organized as follows. In Section 2 we recall some theoretical results: firstly, we introduce the radial basis functions referring to existence and uniqueness of RBF interpolants, then we give a general description of the partition of unity method, which makes use of local RBF approximants. In Section 3, we present in detail the cell algorithm, which is efficiently implemented and optimized by using a cell-based searching procedure. Then, in Section 4 complexity of this algorithm is analyzed as well. Finally, Section 5 shows numerical results concerning efficiency and accuracy of the cell-based partition algorithm. 

\mascotsection{Local interpolation scheme}
\mascotsubsection{Radial basis functions}
A suitable approach to solving the scattered data interpolation problem is to make the assumption that the interpolating function ${\cal R}$ is expressed as a linear combination of \textsl{radial basis functions} $\phi:[0,\infty) \rightarrow \RR$, i.e.,
\begin{eqnarray}
\label{basis}
	{\cal R}(\textbf{x})=\sum_{j=1}^n c_j \phi(||\textbf{x}-\textbf{x}_j||_2), \hspace{1.cm} \textbf{x} \in \RR^2,
\end{eqnarray}
where $||\cdot||_2$ is the Euclidean distance, and ${\cal R}$ satisfies the interpolation conditions \eqref{intcond}.

Thus, solving the interpolation problem under this assumption leads to a system of linear equations of the form
\begin{equation}
	A\textbf{c}=\textbf{f}, \nonumber
\end{equation}
where the entries of the interpolation matrix $A$ are given by 
\begin{eqnarray} \label{coeffmat}
a_{ij}=\phi(||\textbf{x}_i-\textbf{x}_j||_2), \hspace{1 cm} i,j=1,2,\ldots,n,
\end{eqnarray}
$\textbf{c}=\left[c_1,c_2,\ldots,c_n\right]^T$, and $\textbf{f}=\left[f_1,f_2,\ldots,f_n\right]^T$. Moreover, we know that the interpolation problem is \textsl{well-posed}, that is a solution to a problem exists and is unique, if and only if the matrix $A$ is non-singular. A sufficient condition to have non-singularity is that the corresponding matrix is positive definite (see, e.g., \cite{Fasshauer07}).

Now, we remind that a real-valued continuous even function $\phi$ is called positive definite on $\RR^2$ if
\begin{eqnarray} \label{qf}
\sum_{i=1}^{n}\sum_{j=1}^{n} c_i {c_j} \phi\left(\textbf{x}_i - \textbf{x}_j\right) \geq 0
\end{eqnarray} 
for any $n$ pairwise distinct nodes $\textbf{x}_1,\textbf{x}_2,\ldots,\textbf{x}_n \in \RR^2$, and $\textbf{c}=[c_1,c_2,\ldots,$ $c_n]^T\in \RR^n$. The function $\phi$ is called strictly positive definite on $\RR^2$ if the quadratic form (\ref{qf}) is zero only for $\textbf{c}\equiv \textbf{0}$.

Therefore, if $\phi$ is strictly positive definite, the interpolant (\ref{basis}) is unique, since the corresponding interpolation matrix (\ref{coeffmat}) is positive definite and hence non-singular.

Some of the most popular strictly positive definite RBFs are
\begin{equation}
\left.
\begin{array}{rclll}
\phi_G(r) & = & {\rm e}^{-\alpha^2 r^2}, &                      & \hspace{0.5cm} {\rm (Gaussian)} \nonumber \\
\phi_{W2}(r) & = & \displaystyle{\left(1-cr\right)_+^4\left(4cr+1\right)},&   & \hspace{0.5cm} \mbox{{\rm (Wendland's $C^2$ function)}} \nonumber \\
\end{array}
\right.
\end{equation}
where $\alpha, c \in \RR^+$ are the \textsl{shape parameters}, $r=||\textbf{x}-\textbf{x}_i||_2$, and $(\cdot)_+$ denotes the truncated power function. Note that Gaussian (G) is a globally supported function of infinite smoothness, whereas Wendland's function (W2) is compactly supported one with degree of smoothness 2. For further details, see \cite{Wendland05}.

\mascotsubsection{Partition of unity interpolant}
The partition of unity method was firstly suggested in \cite{Babuska97, Melenk96} in the context of meshfree Galerkin methods for the solution of partial differential equations (PDEs), but now it is also commonly used in the field of approximation theory (see \cite{Wendland05}). In particular, this approach enables us to decompose a large problem into many small problems, and at the same time ensures that the accuracy obtained for the local fits is carried over to the global one. In fact, the partition of unity method we here consider can be thought as a modified Shepard's method with higher-order data, where local approximations ${\cal R}_j$ are RBFs. Similar local approaches involving Shepard's type methods were considered in \cite{Allasia11, Lazzaro02, Renka88a, Renka88b, Thacker10}.

Thus, the partition of unity method consists in partitioning the open and bounded domain $\Omega \subseteq \RR^2$ into $d$  subdomains $\Omega_j$ such that $\Omega \subseteq \bigcup_{j=1}^{d} \Omega_j$ with some mild overlap among the subdomains. At first, we choose a partition of unity, i.e. a family of compactly supported, non-negative, continuous functions $W_j$ with $\text{supp}(W_j) \subseteq \Omega_j$ such that 
\begin{equation}
\sum_{j=1}^{d} W_j(\textbf{x}) = 1, \hspace{1. cm} \textbf{x} \in \Omega.   \nonumber
\end{equation}
Then, we can consider the global interpolant
\begin{eqnarray}
\label{pui}
	{\cal I}(\textbf{x})= \sum_{j=1}^{d} {\cal R}_j(\textbf{x}) W_j(\textbf{x}),
\end{eqnarray}
where the local radial basis function
\begin{equation}
	{\cal R}_j(\textbf{x})=\sum_{k=1}^{m_j} c_k \phi(||\textbf{x}-\textbf{x}_k||_2) \nonumber
\end{equation}
is obtained by solving a local interpolation problem, which is constructed using the $m_j$ nodes belonging to each subdomain $\Omega_j$. Note that if the local approximants satisfy the interpolation conditions at node $\textbf{x}_i$, i.e.  ${\cal R}_j(\textbf{x}_i)=f(\textbf{x}_i)$, then the global approximant also interpolates at this node, i.e. ${\cal I}(\textbf{x}_i)=f(\textbf{x}_i)$, for $i=1,2,\ldots,n$. 

In order to be able to formulate error bounds we need some technical conditions. Then, we require the partition of unity functions $W_j$ to be \textsl{k-stable}, i.e. each $W_j \in C^k(\RR^2)$, $j=1,2,\ldots,d$, and for every multi-index $\beta \in \NN_0^m$ with $|\beta| \leq k$ there exists a constant $C_{\beta} > 0$ such that
\begin{equation}
	||D^{\beta}W_j||_{L_{\infty}(\Omega_j)}\leq C_{\beta}/\delta_j^{|\beta|}, \nonumber
\end{equation}
where $\delta_j$ = diam($\Omega_j$).

In accordance with the statements in \cite{Wendland02} we require some additional regularity assumptions on the \textsl{covering} $\{\Omega_j\}_{j=1}^{d}$. Therefore, setting ${\cal X}_n=\{\textbf{x}_i, i=1,2,$ $\ldots,n\} \subseteq \Omega$, an open and bounded covering $\{\Omega_j\}_{j=1}^{d}$ is called \textsl{regular} for $(\Omega,{\cal X}_n)$ if the following properties are satisfied:
\begin{itemize}
	\item[(a)] for each $\textbf{x} \in \Omega$, the number of subdomains $\Omega_j$ with $\textbf{x} \in \Omega_j$ is bounded by a global constant $K$;
	\item[(b)] each subdomain $\Omega_j$ satisfies an interior cone condition;
	\item[(c)] the local fill distances $h_{{\cal X}_{j}, \Omega_j}$ are uniformly bounded by the global fill distance $h_{{\cal X}_n, \Omega}$, where ${\cal X}_{j}={\cal X}_n \cap \Omega_j$.
\end{itemize}
Therefore, assuming that:
\begin{itemize}
	\item  $\phi \in C_{\nu}^k(\RR^2)$ is a strictly positive definite function;
	\item $\{\Omega_j\}_{j=1}^{d}$ is a regular covering for $(\Omega, {\cal X}_n)$; 
	\item $\{W_j\}_{j=1}^{d}$ is $k$-stable for $\{\Omega_j\}_{j=1}^{d}$;
\end{itemize}
we have the following convergence result (see, e.g., \cite{Fasshauer07, Wendland05}), i.e., the error between $f \in {\cal N}_{\phi}(\Omega)$, where ${\cal N}_{\phi}$ is the native space of $\phi$, and its partition of unity interpolant (\ref{pui}) can be bounded by
\begin{equation}
	|D^{\beta}f(\textbf{x}) - D^{\beta}{\cal I}(\textbf{x})| \leq C h_{{\cal X}_n, \Omega}^{(k+\nu)/2 - |\beta|} |f|_{{\cal N}_{\phi}(\Omega)}, \nonumber
\end{equation}
for all $\textbf{x} \in \Omega$ and all $|\beta| \leq k/2$.
 
\mascotsection{Cell algorithm}
In this section we present an efficient algorithm for bivariate interpolation of scattered data sets lying on the domain $\Omega = [0,1]^2 \subset \RR^2$, which is based on the partition of unity method for constructing a global interpolant by blending RBFs as local approximants and using locally supported weight functions. It is efficiently implemented and optimized by connecting the interpolation method itself with an effective cell-based searching procedure. In particular, the implementation of this algorithm is based on the construction of a \textsl{cell structure}, which is obtained by partitioning the domain $\Omega$ in square cells, whose sizes strictly depend on the dimension of its subdomains. Such approach leads to important improvements in the searching processes of the nearest neighbour points compared to the searching techniques presented in \cite{Allasia11, Cavoretto10a, Cavoretto12a}.

\mascotsubsection{Input and output}

\noindent INPUT: 
\begin{itemize}
	\item ${\cal X}_n=\{(x_i,y_i), i=1,2,\ldots,n\}$, set of nodes; 
	\item ${\cal F}_n=\{f_i, i=1,2,\ldots,n\}$, set of data values;  
	\item ${\cal C}_d=\{(\bar{x}_i, \bar{y}_i), i=1,2,\ldots,d\}$, set of subdomain points (centres); 
	\item ${\cal E}_s=\{(\tilde{x}_i,\tilde{y}_i), i=1,2,\ldots,s\}$, set of evaluation points.
\end{itemize}
\vskip3pt
\noindent OUTPUT: 
\begin{itemize}
	\item ${\cal A}_s=\{{\cal I}(\tilde{x}_i,\tilde{y}_i), i=1,2,\ldots,s\}$, set of approximated values.
\end{itemize}

\mascotsubsection{Data partition phase}

\noindent {\tt Stage 1.} The set ${\cal X}_n$ of nodes and the set ${\cal E}_s$ of evaluation points are ordered with respect to a common direction (e.g. the $y$-axis), by applying a \textsl{quicksort$_y$ procedure}.
\vskip 3pt
\noindent {\tt Stage 2.} For each subdomain point $(\bar{x}_i,\bar{y}_i)$, $i=1,2,\ldots,d$, a local circular subdomain is constructed, whose half-size (the radius) depends on the subdomain number $d$, that is
\begin{eqnarray}
\label{delta}
 \delta_{subdom} = \sqrt{\frac{2}{d}}.
 \end{eqnarray}
This value is suitably chosen, supposing to have a nearly uniform node distribution and assuming that the ratio $n/d \approx 4$.
\vskip 3pt 
\noindent {\tt Stage 3.} A double structure of crossed strips is constructed as follows:
\begin{enumerate}
	\item[i)] a first family of $q$ strips, parallel to the $x$-axis, is considered taking
\begin{eqnarray} \label{q_str}
	q = \left\lceil \frac{1}{\delta_{subdom}} \right\rceil,
\end{eqnarray}
and a \textsl{quicksort$_x$ procedure} is applied to order the nodes belonging to each strip;
	\item[ii)] a second family of $q$ strips, parallel to the $y$-axis, is considered.
\end{enumerate} 
Note that each of the two strip structures are ordered and numbered from 1 to $q$; moreover, the choice in (\ref{q_str}) follows directly from the side length of the domain $\Omega$ (unit square), that here is $1$, and the subdomain radius $\delta_{subdom}$.
\vskip 3pt 
\noindent {\tt Stage 4.} The domain (unit square) is partitioned by a \textsl{cell-based structure} consisted of $q^2$ square cells, whose length of the sides is given by $\delta_{cell} \equiv \delta_{subdom}$. Then, the following structure is considered: 
\begin{itemize}
	\item the sets ${\cal X}_n$, ${\cal C}_d$ and ${\cal E}_s$ are partitioned by the cell structure into $q^2$ subsets ${\cal X}_{n_k}$, ${\cal C}_{d_k}$ and ${\cal E}_{p_k}$, $k=1,2,\ldots,q^2$,
\end{itemize}
where $n_k$, $d_k$ and $p_k$ are the number of points in the $k$-th cell.

\mascotsubsection{Localization phase}

\noindent {\tt Stage 5.} In order to identify the cells to be examined in the searching procedure, we consider two steps as follows: 
\begin{enumerate}
\item[(A)] since $\delta_{cell} \equiv \delta_{subdom}$, the ratio between these quantities is denoted by $i^* = \delta_{subdom}/\delta_{cell}=1$. So the number $j^* = (2i^*+1)^2$ of cells to be examined for each node is $9$.
\item[(B)] for each cell $k=[v,w]$, $u,v=1,2,\ldots,q$, a cell-based searching procedure is considered, examining the points from the cell $[v-i^*,w-i^*]$ to the cell $[v+i^*,w+i^*]$. Note that if $v-i^* < 1$ and/or $w-i^* < 1$, or $v+i^* > q$ and/or $w+i^* > q$, then we set $v-i^* = 1$ and/or $w-i^* = 1$, and $v+i^* = q$ and/or $w+i^* = q$. 
\end{enumerate}

Then, after defining which and how many cells are to be examined, the cell-based searching procedure is applied:
\begin{itemize}
	\item for each subdomain point of ${\cal C}_{d_k}$, $k=1,2,\ldots,q^2$, to determine all nodes belonging to a subdomain. The number of nodes of the subdomain centred at $(\bar{x}_i,\bar{y}_i)$ is counted and stored in $m_i$, $i=1,2,\ldots,d$;
	\item for each evaluation point of ${\cal E}_{p_k}$, $k=1,2,\ldots,q^2$, in order to find all those belonging to a subdomain of centre $(\bar{x}_i,\bar{y}_i)$ and radius $\delta_{subdom}$. The number of subdomains containing the $i$-th evaluation point is counted and stored in $r_i$, $i=1,2,\ldots,s$. 
\end{itemize}

\mascotsubsection{Evaluation phase}

\noindent {\tt Stage 7.} A local approximant $R_j(x,y)$ and a weight function $W_j(x,y)$, $j=1,2,\ldots,d$, is found for each evaluation point.
\vskip 3pt 
\noindent {\tt Stage 8.} Applying the global fit (\ref{pui}), the surface can be approximated at any evaluation point $(\tilde{x},\tilde{y}) \in {\cal E}_s$.

\mascotsection{Complexity}
The partition of unity algorithm involves the use of the standard sorting routine quicksort, which requires on average a time complexity ${\cal O}(M\log M)$, where $M$ is the number of nodes to be sorted. Specifically, we have a data partition phase consisting of building the data structure, where the computational cost is: 
\begin{itemize}
	\item ${\cal O}(n\log n)$ for sorting all $n$ nodes;
	\item ${\cal O}(s\log s)$ for sorting all $s$ evaluation points.
\end{itemize}
Moreover, in order to compute the local RBF interpolants, we have to solve $d$ linear systems and the cost is: 
\begin{itemize}
	\item  ${\cal O}(m_i^3)$, $i=1,2,\ldots,d$, where $m_i$ is the number of nodes in the $i$-th subdomain. 
\end{itemize}
Then, for the $k$-th evaluation point of ${\cal E}_s$ the cost is:
\begin{itemize}
	\item  $r_k \cdot {\cal O}(m_i)$, $i=1,2,\ldots,d$, $k=1,2,\ldots,s$.
\end{itemize}
Finally, the algorithm requires $3n$, $3d$ and $3s$ storage requirements for the data, and $m_i$, $i=1,2,\ldots, d$, locations for the coefficients of each local RBF interpolant.

\mascotsection{Numerical results}
In this section we present some tests to verify performance and effectiveness of the cell-based partition algorithm on scattered data sets. The code is implemented in C/C++ language, while numerical results are carried out on a Intel Core 2 Duo Computer (2.1 GHz). In the experiments we consider a node distribution with $n=4225,16641,66049$ uniformly random Halton nodes generated by using the program given in \cite{Wong97}. The partition of unity algorithm is run considering $d=1024,4096,16384$ subdomain points and $s = 33 \times 33$ evaluation (or grid) points, which are contained in the unit square $\Omega = [0,1]^2$

The performance of the interpolation algorithm is verified taking the data values by Franke's test function
\vskip 0.3cm
$\displaystyle{f(x,y)=\frac{3}{4}\exp\left[-\frac{(9x-2)^2+(9y-2)^2}{4}\right]+\frac{3}{4}\exp\left[-\frac{(9x+1)^2}{49}-\frac{9y+1}{10}\right]}$

$\hspace{1.cm}\displaystyle{+\frac{1}{2} \exp\left[-\frac{(9x-7)^2+(9y-3)^2}{4}\right]-\frac{1}{5} \exp\left[-(9x-4)^2-(9y-7)^2\right]}$.
\vskip 0.3cm

Moreover, since we are concerned to point out the effectiveness of the proposed algorithm, in Table \ref{time_comparison} we compare CPU times (in seconds) obtained by running the cell algorithm described in Section 3, and the strip algorithm proposed in \cite{Cavoretto12c}. This comparison highlights the high efficiency of the cell algorithm, which gives us a considerable saving of time. 

\begin{table}[ht!]
		\begin{center}
			\begin{tabular}{ccccc} 
			 \hline
			 \rule[0mm]{0mm}{3ex}
			  	$n$    & $d$ & \textbf{$t_{cell}$} & $t_{strip}$ \\
				\hline
				\rule[0mm]{0mm}{3ex}
				$4225$    & $1024$   & $\textbf{0.3}$  & $0.4$        \\
				\rule[0mm]{0mm}{3ex}
				$16641$   & $4096$   & $\textbf{0.8}$  & $1.3$        \\
			  \rule[0mm]{0mm}{3ex}
			  $66049$   & $16384$  & $\textbf{2.6}$  & $6.5$       \\
			  \hline 
			\end{tabular}
		\end{center}
			\caption{CPU times (in seconds) obtained by running the cell algorithm ($t_{cell}$) and the strip algorithm ($t_{strip}$).}
			\label{time_comparison}
	\end{table}

Then, in order to test accuracy of the local algorithm, in Table \ref{tab1_errors} we report the Root Mean Square Errors (RMSEs), i.e.
\begin{equation}
	RMSE = \sqrt{\frac{1}{s}\sum_{i=1}^{s} |f(\textbf{x}_i) - {\cal I}(\textbf{x}_i)|^2}. \nonumber
\end{equation}
The error computation is achieved by considering both globally and locally supported RBFs for suitable values of the shape parameters, i.e., $\alpha^2 = 50$ for $\phi_G$, and $c = 1$ for $\phi_{W2}$. We observe that the local scheme turns out to be accurate, even if we do not consider the optimal values for the shape parameters, namely those values for which we get the best possible results. However, these choices give a good compromise among accuracy and stability. 

\begin{table}[ht!]
		\begin{center}
			\begin{tabular}{cccc} 
			 \hline
			 \rule[0mm]{0mm}{3ex}
			  	$n$     & $4225$ & $16641$ & $66049$ \\
				\hline
				\rule[0mm]{0mm}{3ex}
				  $\phi_G$  & $3.0045{\rm E}-4$  & $2.8214{\rm E}-5$ & $1.5200{\rm E}-6$      \\
			  \rule[0mm]{0mm}{3ex}
				  $\phi_{W2}$  & $2.2145{\rm E}-4$  & $5.3127{\rm E}-5$ & $9.3027{\rm E}-6$      \\
			  \hline 
			\end{tabular}
		\end{center}
			\caption{RMSEs with $\alpha^2 = 50$ and $c = 1$ for Franke's function.}
			\label{tab1_errors}
	\end{table}

\begin{figure}[ht!]
\begin{center}
\begin{minipage}{60mm}
\includegraphics[width=6.0cm]{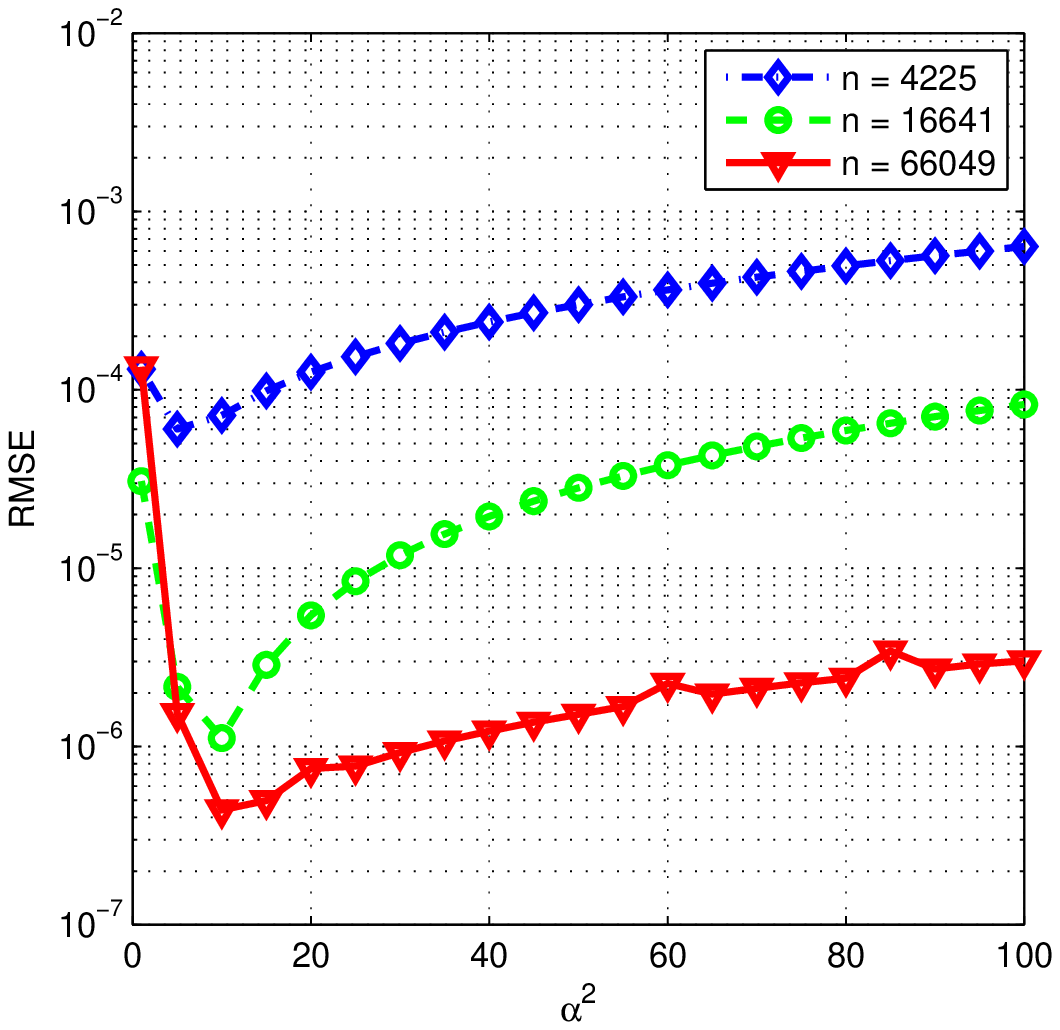}
\centerline{$\phi_G$}
\end{minipage}
\hfil
\begin{minipage}{60mm}
\includegraphics[width=6.0cm]{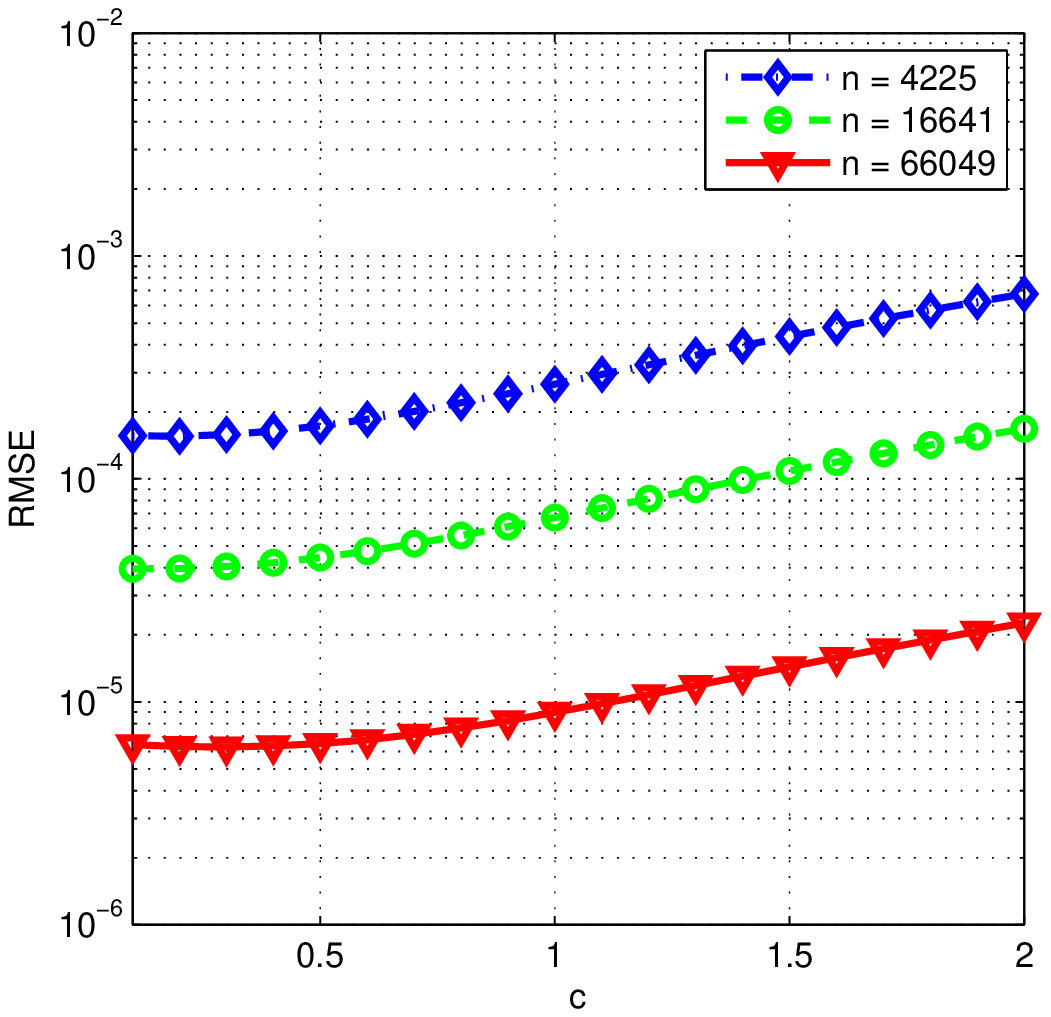}
\centerline{$\phi_{W2}$}
\end{minipage}
\end{center}
\caption{RMSEs obtained by varying $\alpha^2$ and $c$ for Franke's function.}
\label{shape_f1}
\end{figure}

Then, in Figure \ref{shape_f1} we plot the behavior of the RMSEs by varying values of the shape parameters for Franke's function. These graphs (and other ones obtained considering different test functions we omit for brievity) point out that, if an optimal search of the shape parameters was performed, in some cases the results of accuracy reported in this section could be improved of one or even two orders of magnitude. Note that each evaluation is carried out by choosing equispaced values of the shape parameter with $\alpha^2 \in [1, 100]$, and $ c \in [0.1, 2]$.

By analyzing numerical tests and the related pictures, we observe that Wendland's function $\phi_{W2}$ has greater stability than Gaussian $\phi_G$ and good accuracy. However, these graphs give an idea on stability and enable us to choose \lq\lq sure\rq\rq~values for the shape parameters. These tests confirm theoretical results and suggest to use basis functions with a moderate order of smoothness, thus avoiding the well-known ill-conditioning problems of infinitely smooth RBFs, in particular if we are dealing with a very large number of nodes. 

\end{document}